\documentclass{article}
\usepackage{stmaryrd}
\usepackage{mathrsfs}
\usepackage{amsmath}
\usepackage{amssymb}
\usepackage{amsfonts}
\usepackage{pstricks}
\usepackage{graphicx, color, epsfig, enumerate}

\newcommand{\om}{\omega}
\newcommand{\ga}{\gamma}

\newcommand{\la}{\lambda}
\newcommand{\q}{\tilde}
\newcommand{\al}{\alpha}

\newcommand{\f}{\frac}

\newtheorem{theorem}{Theorem}[section]
\newtheorem{proposition}[theorem]{Proposition}

\numberwithin{equation}{section}

\begin{document}
\title{{A Novel Approach to Elastodynamics:}\\
{I. The Two-Dimensional Case}}
\author{
{A. S. Fokas${}^a$\footnote{T.Fokas@damtp.cam.ac.uk}\  \ and
D. Yang${}^{a,}{}^{b}$\footnote{yangd04@mails.tsinghua.edu.cn}}\\
\\
{\small ${}^a$ Department of Applied Mathematics and Theoretical Physics,}\\
{\small University of Cambridge, Cambridge CB3 0WA, UK}\\
\\
{\small ${}^b$ Department of Mathematical Sciences, Tsinghua
University,}\\
{\small Beijing 100084, P. R. China}\\
}
\maketitle

\begin{abstract}
We introduce a new approach to constructing analytic solutions of
the linear PDEs describing elastodynamics. This approach is
illustrated for the case of a homogeneous isotropic half-plane body
satisfying arbitrary initial conditions and Lamb's boundary
conditions. A particular case of this problem, namely the case of
homogeneous initial conditions, was first solved by Lamb using the
Fourier-Laplace transform. The solution of the general problem can
also be expressed in terms of the Fourier transform, but this
representation involves transforms of \textit{unknown} boundary
values. This necessitates the formulation and solution of a
cumbersome auxiliary problem, which expresses the unknown boundary
values in terms of the Laplace transform of the given boundary
data. The new approach, which is applicable to arbitrary
initial and boundary conditions, bypasses the above auxiliary
problem and expresses the solutions directly in terms of the given
initial and boundary data.
\end{abstract}
\noindent {\small{\sc Keywords}: elastodynamics, half space, initial
boundary value problem, Lamb's problem, global relation.}

\section{Introduction}
The problem considered in this paper has a long and illustrious history,
which begins with the classic works of Sir Horace Lamb in 1904
\cite{HL, JM, KLJ}. In \cite{HL}, Lamb treated four basic
problems, the so called Lamb's problems, including both the two and the three dimensional cases. 
These problems have been studied by several authors, see for example
\cite{point_load0, point_load05, point_load1, point_load2,
point_load4, point_load5, line_load1, line_load16, line_load2}. For
the axisymmetric point-load problem, Pekeris \cite{point_load0}
obtained an exact expression for the surface motion in the case that
the point load is given by a Heaviside function. This result was
further extended by Schiel and Prot\'azio \cite{point_load1}. For
the three dimensional cases, Payton, using the Betti-Rayleigh
reciprocal theorem \cite{point_load05}, obtained the displacement
and the surface motion. Gakenheimer, in his PhD thesis
\cite{point_load2}, made a systematic investigation of the
travelling normal point load, by using the Cagniard-deHoop method.
Bakkeraq et al. in \cite{point_load5} revisited this problem and
presented a straightforward implementation of the Cagniard-deHoop
method. Several related problems of contact mechanics and wave
propagation can be found in the books of Miklowitz \cite{JM}, of
Johnson \cite{KLJ}, and of Payton \cite{Payton}.

The extensive study of Lamb's problems is perhaps due to the fact
that they arise in a variety of physical situations, as well as to
their interesting mathematical structure.

In this paper, we study two-dimensional problems under the usual
assumption of plane strain. We consider arbitrary initial conditions
and general stress boundary conditions, including normal line load,
tangential line load, and mixed line load. We refer to the latter
three stress conditions as Lamb's boundary conditions.

Most studies of Lamb's problems are based on the Helmholtz
decomposition and on the use of the Laplace transform in time. In
particular, Lamb was the first to use the above approach, as well as
a certain superposition which is in fact equivalent to the use of
the Fourier transform. Miklowitz made use of the Cagniard-deHoop
method \cite{JM}, which itself is based on the Laplace transform.
Payton made use of Green's functions and of the Laplace transform.

Helmholtz decomposition has the advantage of decomposing the P-waves
and S-waves. However, it has the disadvantage of introducing higher
order derivatives to the boundary conditions. The use of Laplace
transform in time, essentially restricts the problem to the case of
homogeneous initial conditions.

Our new approach is based on the unified approach for solving linear and
integrable nonlinear PDEs introduced by one of the authors in \cite{ASF1997}.
A crucial role in this method is played by the so-called \textit{global
relations}, which are algebraic relations coupling appropriate transforms of the
unknown boundary values with transforms of the given data.
For linear PDEs this method uses three novel steps \cite{ASF}-\cite{SSF}: 1. Derive a representation for the solution
in terms of an integral involving a contour in the complex Fourier
plane. This representation is not yet effective, because in addition to
transforms of the given initial and boundary conditions, it also
contains transforms of unknown boundary values. 2. Analyse certain
transformations in the complex Fourier plane which leave invariant the
transforms of the unknown boundary values. 3. Eliminate the transforms
of the unknown boundary values, by combining the results of the first two
steps and by employing Cauchy's theorem (or more precisely Jordan's lemma).

This paper is organised as follows. In section 2, we recall the
governing equations of elastodynamics and derive the global relations. In section 3 we
implement step 1 and in section 4 we implement steps 2 and 3. In section 5 we use
the general representations derived in section 4 in order to analyse
Lamb's problems. These results are further discussed in section 6.

For certain complicated boundary value problems it seems that it is not possible
to eliminate from the integral representation of the solution the
transforms of the unknown boundary values. However, for some of these
problems, by using the global relations, one can derive expressions for
the Laplace transforms of the boundary values in terms of the given data \cite{Ashton_ASF}.
A summary of this \textit{less} effective approach is presented in
the Appendix. It is interesting that for the \textit{particular case of zero
initial conditions}, the formulae presented in the Appendix reduce to the
formulae first derived in the classical works of Lamb.
\section{Governing Equations and Global Relations}
The transient problem for two dimensional elastodynamics in the half plane with
Lamb's boundary conditions is defined as follows: Let
$u=u(x,y,t),$ $v=v(x,y,t),$ denote the displacements of a
homogeneous isotropic half space body. The governing equations of
motion without external body forces, are the Lam\'e-Navier
equations:
\begin{subequations}\label{PLN}
\begin{align}
\label{PLN1}
(\la+2\mu)u_{xx}+(\la+\mu)v_{xy}+\mu u_{yy}-\rho u_{tt}=0,\\
\label{PLN2}
\mu v_{xx}+(\la+\mu)u_{xy}+(\la+2\mu)v_{yy}-\rho v_{tt}=0,
\end{align}
\end{subequations}
\begin{flushright}
 $-\infty<x<\infty,~0<y<\infty,~t>0,$
\end{flushright}
where $\la,\mu$ are the Lam\'e constants and $\rho$ denotes the density of
the material, which without loss of generality is normalized to unity, i.e., $\rho\equiv1.$
Let the initial conditions be denoted by
\begin{subequations}\label{ic}
\begin{align}
\label{ic1}u(x,y,0)&=u_0(x,y),~u_t(x,y,0)=u_1(x,y),\\
\label{ic2}v(x,y,0)&=v_0(x,y),~v_t(x,y,0)=v_1(x,y),
\end{align}
\end{subequations}
\begin{flushright}
 $-\infty<x<\infty,~0<y<\infty.$
\end{flushright}

Let the stress boundary conditions be denoted by
\begin{subequations}\label{bv}
\begin{align}
\label{bv1}(u_y+v_x)(x,0,t)&=g_1(x,t),~~-\infty<x<\infty,~t>0,\\
\label{bv2}
\Big(v_y+\f{\la}{\la+2\mu}u_x\Big)(x,0,t)&=g_2(x,t),~~-\infty<x<\infty,~t>0.
\end{align}
\end{subequations}
For a tangential line load, the functions $g_1$ and $g_2$ are given by
\begin{equation}
g_1(x,t)=\sigma_0\delta(x)X(t)/\mu,~g_2(x,t)=0;
\end{equation}
for a normal line load,
\begin{equation}
g_1(x,t)=0,~g_2(x,t)=\sigma_0\delta(x)Y(t)/(\la+\mu);
\end{equation}
for a moving normal line load with a constant speed $C$,
\begin{equation}
g_1(x,t)=0,~g_2(x,t)=\sigma_0\delta(x-Ct)/(\la+\mu).
\end{equation}
Here $\sigma_0$ is a constant, $\delta(x)$ denotes the Dirac-$\delta$ function, $X(t)$ and $Y(t)$ are functions which depend only on $t$.
\paragraph{Notations}
Hat, ``$\wedge$'', will denote the two-dimensional Fourier transform with respect to
$x$ and $y$, whereas tilde, ``$\sim$'', will denote the Fourier
transform with respect to $x$. In particular,
\begin{subequations}\label{Fourier}
\begin{align}
&\hat{u}(k,l,t):=\int_{-\infty}^{\infty}dx\int_0^{\infty}dy~
e^{-ikx-ily}u(x,y,t),\\
&\hat{v}(k,l,t):=\int_{-\infty}^{\infty}dx\int_0^{\infty}dy~
e^{-ikx-ily}v(x,y,t),\\
&\hat{u}_j(k,l):=\int_{-\infty}^{\infty}dx\int_0^{\infty}dy~
e^{-ikx-ily}u_j(x,y),\\
&\hat{v}_j(k,l):=\int_{-\infty}^{\infty}dx\int_0^{\infty}dy~
e^{-ikx-ily}v_j(x,y),
\end{align}
\begin{flushright}$k\in\mathbb{R},~l\in\mathbb{C}^-,~t>0, j=0,1,$\end{flushright}
\end{subequations}
where $\mathbb{C}^-$ denotes the lower half complex $l$-plane.

Furthermore,
\begin{subequations}\label{xF}
\begin{align}
\tilde{u}(k,t)&:=\int_{-\infty}^{\infty}dx~
e^{-ikx}u(x,0,t),\\
\tilde{v}(k,t)&:=\int_{-\infty}^{\infty}dx~
e^{-ikx}v(x,0,t), ~~ k\in\mathbb{R},t>0.
\end{align}
\end{subequations}
\begin{equation}\label{gtilde}
\q{g}_j(k,t):=\int_{-\infty}^{\infty}dx~e^{-ikx} g_j(x,t),~~k\in\mathbb{R},~t>0,~j=1,2.
\end{equation}

We emphasise that since $x\in \mathbb{R}$, the $x-$Fourier transform
is well-defined only for $k\in\mathbb{R};$ on the other hand, since
$0<y<\infty$,  the $y-$Fourier transform is well-defined for $l$ 
in the lower half complex $l-$plane.

By applying the two-dimensional Fourier transform to equations \eqref{PLN} and by using in the resulting equations
the boundary conditions \eqref{bv}, we obtain the following equations:
\begin{subequations}\label{pre}
\begin{align}\label{pre_g1}
&-(\la+2\mu) k^2\hat{u}+(\la+\mu)ik(il\hat{v}-\q{v})+\mu\big[-l^2\hat{u}-il\q{u}-(\q{g}_1-ik\q{v})\big]=\hat{u}_{tt},\\
\label{pre_g2}
\begin{split}
&-\mu k^2\hat{v}+(\la+\mu)ik(il\hat{u}-\q{u})+(\la+2\mu)\Big[-l^2\hat{v}-il\q{v}
-\Big(\q{g}_2-\f{\la}{\la+2\mu}ik\q{u}\Big)\Big]\\
&~~~~~~~~~~~~~~~~~~~~~~~~~~~~~~~~~~~~~~~~~~~~~~~~~~~~~~~~~~~~~~=\hat{v}_{tt},~~~~~~ k\in\mathbb{R},~l\in\mathbb{C}.
\end{split}
\end{align}
\end{subequations}
Introducing the notations
\begin{subequations}
\begin{align}
 P(k,l,t)=k\hat{u}(k,l,t)+l\hat{v}(k,l,t),\\
 Q(k,l,t)=l\hat{u}(k,l,t)-k\hat{v}(k,l,t),
\end{align}
\end{subequations}
equations \eqref{pre} become
\begin{subequations}\label{ODEs}
\begin{align}
&P_{tt}+(\la+2\mu)(k^2+l^2)P=F_P,\\
&Q_{tt}+\mu(k^2+l^2)Q=F_Q,~~~k\in\mathbb{R},~l\in\mathbb{C},
\end{align}
\end{subequations}
where the functions $F_P(k,l,t)$ and $F_Q(k,l,t)$ are defined as follows:
\begin{subequations}\label{Bpq}
\begin{align}
\label{Bp}
\begin{split} F_P(k,l,t)=&-[\mu k\q{g}_1(k,t)+i\la
k^2\q{v}(k,t)]-l[(\la+2\mu)\q{g}_2(k,t)\\
&+2i\mu k\q{u}(k,t)]-il^2(\la+2\mu)\q{v}(k,t),
\end{split}\\
\begin{split}
\label{Bq}
F_Q(k,l,t)=&[k(\la+2\mu)\q{g}_2(k,t)+i\mu
k^2\q{u}(k,t)]-l[\mu\q{g}_1(k,t)\\
&-2i\mu k\q{v}(k,t)]-i\mu l^2\q{u}(k,t).
\end{split}
\end{align}
\end{subequations}
Solving equations \eqref{ODEs} for $\{P,Q\}$, we obtain the following expressions:
\begin{subequations}\label{global_r}
 \begin{align}
\label{global_r1}
\begin{split}
k\hat{u}+l\hat{v}=&\f{1}{-2i\om_1}\Big(e^{-i\om_1t}\int_0^t
e^{i\om_1s}F_P(k,l,s)ds-e^{i\om_1t}\int_0^t e^{-i\om_1s}F_P(k,l,s)ds\Big)\\
&+\Big(\f{1}{2}P_0+\f{i}{2}\f{P_1}{\om_1}\Big)e^{-i\om_1t}+\Big(\f{1}{2}P_0-\f{i}{2}\f{P_1}{\om_1}\Big)e^{i\om_1t},
\end{split}\\
\label{global_r2}
\begin{split}
l\hat{u}-k\hat{v}=&\f{1}{-2i\om_2}\Big(e^{-i\om_2t}\int_0^t
e^{i\om_2s}F_Q(k,l,s)ds-e^{i\om_2t}\int_0^t e^{-i\om_2s}F_Q(k,l,s)ds\Big)\\
&+\Big(\f{1}{2}Q_0+\f{i}{2}\f{Q_1}{\om_2}\Big)e^{-i\om_2t}+\Big(\f{1}{2}Q_0-\f{i}{2}\f{Q_1}{\om_2}\Big)e^{i\om_2t},
\end{split}
 \end{align}
\end{subequations} \begin{flushright} $k\in\mathbb{R},~l\in\mathbb{C},$\end{flushright}
where the dispersion relations $\om_1$ and $\om_2$ are given by
\begin{equation}\label{dispersion}
\om_1^2=(\la+2\mu)(k^2+l^2),~~\om_2^2=\mu(k^2+l^2),~~k\in\mathbb{R},~l\in\mathbb{C},
\end{equation}
and the \textit{known} functions $P_0(k,l),P_1(k,l),Q_0(k,l),Q_1(k,l)$ are given in terms of the
initial conditions by
\begin{subequations}\label{initial}
 \begin{align}
P_0=k \hat{u}_0+l\hat{v}_0,~P_1=k\hat{u}_1+l\hat{v}_1,\\
Q_0=l \hat{u}_0-k\hat{v}_0,~P_1=l\hat{u}_1-k\hat{v}_1,
 \end{align}
\end{subequations}
with $\hat{u}_j,~\hat{v}_j,~j=0,1,$ defined in equations \eqref{Fourier}.

In what follows, we take
\begin{equation}\label{dispersion_choice}
\om_1=\sqrt{\la+2\mu}(k^2+l^2)^{\f{1}{2}},~~\om_2=\sqrt{\mu}(k^2+l^2)^{\f{1}{2}}.
\end{equation}
The function $(k^2+l^2)^{\f{1}{2}}$ has the branch points $\pm
ik$ in the complex $l-$plane; we connect these two branch points by
a branch cut and we fix a branch in the cut plane by the requirement
that,
$$(k^2+l^2)^{\f{1}{2}}\sim~ l+O\Big(\f{1}{l}\Big),~~~as~~l\rightarrow\infty.$$

Let $u^{(j)\pm},~v^{(j)\pm},~U^{(j)},~V^{(j)},~j=1,2,$ denote the
following \textit{unknown} functions:
\begin{subequations}\label{BigU}
 \begin{align}
&u^{(j)\pm}(k,l,t)=\int_0^t e^{\pm i\om_j s}\q{u}(k,s)ds,\\
&v^{(j)\pm}(k,l,t)=\int_0^t e^{\pm i\om_j s}\q{v}(k,s)ds,\\
\label{defU1}&U^{(j)}(k,l,t)=\f{1}{2\om_j}\big(e^{-i\om_j t}u^{(j)+}(k,l,t)-e^{i\om_j t}u^{(j)-}(k,l,t)\big),\\
&V^{(j)}(k,l,t)=\f{1}{2\om_j}\big(e^{-i\om_j
t}v^{(j)+}(k,l,t)-e^{i\om_j t}v^{(j)-}(k,l,t)\big),
 \end{align}
\begin{flushright}
$k\in\mathbb{R}, l\in\mathbb{C}, t\geq0, j=1,2.$\end{flushright}
\end{subequations}
Similarly, let $g^{(j)\pm}$, $f^{(j)\pm}$, $G^{(j)}$,
$F^{(j)},~j=1,2,$ denote the following \textit{known} functions:
\begin{subequations}\label{BigG}
 \begin{align}
&g^{(j)\pm}(k,l,t)=\int_0^t e^{\pm i\om_j s}\q{g}_1(k,s)ds,\\
&f^{(j)\pm}(k,l,t)=\int_0^t e^{\pm i\om_j s}\q{g}_2(k,s)ds,\\
&G^{(j)}(k,l,t)=\f{1}{2\om_j}\big(e^{-i\om_j t}g^{(j)+}(k,l,t)-e^{i\om_j t}g^{(j)-}(k,l,t)\big),\\
&F^{(j)}(k,l,t)=\f{1}{2\om_j}\big(e^{-i\om_j
t}f^{(j)+}(k,l,t)-e^{i\om_j t}f^{(j)-}(k,l,t)\big),
\end{align}
\begin{flushright}$k\in\mathbb{R}, l\in\mathbb{C}, t\geq0, j=1,2.$\end{flushright}
\end{subequations}

Using the above notations, equations \eqref{global_r} become
\begin{subequations}\label{global_relations}
 \begin{align}
\label{1} k\hat{u}+l\hat{v}=&(\la k^2+l^2(\la+2\mu))V^{(1)}+ 2\mu k l U^{(1)}+N_P,\\
\label{2} l\hat{u}-k\hat{v}=&-\mu(k^2-l^2)U^{(2)}-2\mu k l V^{(2)}+N_Q,
 \end{align}
\end{subequations}
where the \textit{known} functions $N_P(k,l,t)$ and $N_Q(k,l,t)$ are defined as follows:
\begin{subequations}\label{NPNQ}
 \begin{align}
\begin{split}
N_P(k,l,t)=&-il(\la+2\mu) F^{(1)}(k,l,t) - i\mu k G^{(1)}(k,l,t)\\
+&\Big(\f{1}{2}P_0(k,l)+\f{i}{2}\f{P_1(k,l)}{\om_1}\Big)e^{-i\om_1t}+\Big(\f{1}{2}P_0(k,l)-\f{i}{2}\f{P_1(k,l)}{\om_1}\Big)e^{i\om_1t},\\
\end{split}\\
\begin{split}
N_Q(k,l,t)=&ik(\la+2\mu) F^{(2)}(k,l,t)- i\mu l G^{(2)}(k,l,t)\\
&+\Big(\f{1}{2}Q_0(k,l)+\f{i}{2}\f{Q_1(k,l)}{\om_2}\Big)e^{-i\om_2t}+\Big(\f{1}{2}Q_0(k,l)-\f{i}{2}\f{Q_1(k,l)}{\om_2}\Big)e^{i\om_2t}.\\
\end{split}
\end{align}
\end{subequations}
We will refer to equations \eqref{global_relations} as the \textit{global relations}.
The global relations express the two-dimensional Fourier
transforms of the solution $(u,v)$ in terms of the given initial and
boundary data, as well as in terms of the transforms
$U^{(1)},U^{(2)},V^{(1)},V^{(2)}$ of the unknown boundary values. The important
observation is that equations \eqref{global_relations} are valid for all
values of $l$ in the lower half complex $l-$plane. It turns out that
using this fact it will be possible to eliminate the unknown
transforms.

\section{An Integral Representation Involving the Unknown Transforms}
We first observe that for any fixed $k\in\mathbb{R}$ and any fixed
$t$, $0\leq t< T$, $T>0$, the functions $U^{(j)}$, $V^{(j)}$,
$G^{(j)},F^{(j)},$ $j\in\{1,2\},$ are analytic in the complex
$l-$plane.

Indeed, these functions are all single-valued, thus it only remains
to establish the analyticity in the neighbourhood $l=\pm ik$. The
definition of $U^{(1)}$, i.e. equation \eqref{defU1}, implies that
this function possesses the following expansion for $l$ near $\pm ik$:
\begin{equation}\label{expansion}
U^{(1)}(k,l,t)=i\sum_{n=1}^{\infty}\f{\om_1^{2n-2}}{(2n-1)!}\int_0^t(s-t)^{2n-2}\q{u}(k,s)ds.
\end{equation}
Similar expansions are also valid for $U^{(2)},V^{(1)},V^{(2)},
G^{(1)},G^{(2)},F^{(1)},F^{(2)}.$

Solving equations \eqref{global_relations} we find
\begin{subequations}
 \begin{align}
\label{gwh1}
\begin{split}
\hat{u}=\f{1}{k^2+l^2}\Big\{& (\la k^3+k l^2(\la+2\mu)) V^{(1)}+ 2\mu k^2 l U^{(1)}\\
&-\mu( k^2l-l^3) U^{(2)} -2\mu k l^2V^{(2)}+kN_P+lN_Q \Big\},
\end{split}\\
\label{gwh2}
\begin{split}
\hat{v}=\f{1}{k^2+l^2}\Big\{& (\la k^2 l+ l^3(\la+2\mu))V^{(1)} + 2\mu k l^2 U^{(1)}\\
&+\mu(k^3-k l^2) U^{(2)} + 2\mu k^2 lV^{(2)}+l N_P - k N_Q \Big\}.
\end{split}
 \end{align}
\end{subequations}

We observe that the following important transformation is valid in
the complex $l-$plane:
\begin{equation}\label{superb_1}
\begin{split}
l\rightarrow -l~~&maps~~\om_1~~ to~~ -\om_1,~~~ \om_2~~ to~~ -\om_2,~
and~~\\
&leaves~~ U^{(j)},V^{(j)},~ j=1,2~~ invariant.
\end{split}
\end{equation}
Employing the transformation $l\rightarrow-l$ in equations
\eqref{gwh1} and \eqref{gwh2} and then adding the resulting
equations to \eqref{gwh1} and \eqref{gwh2}, we obtain the following
equations:
\begin{subequations}
 \begin{align}
\label{u}
\begin{split}
&\hat{u}(k,l,t)+\hat{u}(k,-l,t)=\\
&~~~~\f{1}{k^2+l^2}\Big\{ 2(\la k^3+ k l^2(\la+2\mu))V^{(1)}(k,l,t)-4\mu k l^2 V^{(2)}(k,l,t)+\\
&~~~~k(N_P(k,l,t)+N_P(k,-l,t))+l(N_Q(k,l,t)-N_Q(k,-l,t)) \Big\},
\end{split}\\
\label{v}
\begin{split}
&\hat{v}(k,l,t)+\hat{v}(k,-l,t)=\\
&~~~~\f{1}{k^2+l^2}\Big\{ 4\mu k l^2 U^{(1)}(k,l,t)+2\mu(k^3 - k l^2) U^{(2)}(k,l,t)+\\
&~~~~l(N_P(k,l,t)-N_P(k,-l,t))-k(N_Q(k,l,t)+N_Q(k,-l,t)) \Big\}.
\end{split}
 \end{align}
\end{subequations}
Applying the inverse Fourier transform formula to these equations,
we obtain
\begin{subequations}\label{first}
 \begin{align}
\label{first_u}
\begin{split}
&u(x,y,t)=\f{1}{4\pi^2}\int_{-\infty}^{\infty}dk \int_{-\infty}^{\infty} dl~e^{ikx+ily}\hat{u}(k,l,t)\\
=&\f{1}{4\pi^2}\int_{-\infty}^{\infty}dk \int_{-\infty}^{\infty} dl
\f{e^{ikx+ily}}{k^2+l^2} \Big\{ 2(\la k^3+k l^2(\la+2\mu))V^{(1)}(k,l,t)-4\mu k l^2 V^{(2)}(k,l,t)\\
&+k(N_P(k,l,t)+N_P(k,-l,t))+l(N_Q(k,l,t)-N_Q(k,-l,t))\Big\},
\end{split}\\
\label{first_v}
\begin{split}
&v(x,y,t)=\f{1}{4\pi^2}\int_{-\infty}^{\infty}dk \int_{-\infty}^{\infty} dl~ e^{ikx+ily}\hat{v}(k,l,t)\\
=&\f{1}{4\pi^2}\int_{-\infty}^{\infty}dk \int_{-\infty}^{\infty} dl
\f{e^{ikx+ily}}{k^2+l^2}\Big\{4\mu k l^2 U^{(1)}(k,l,t)+2(\mu k^3-\mu k l^2) U^{(2)}(k,l,t)\\
&+l(N_P(k,l,t)-N_P(k,-l,t))-k(N_Q(k,l,t)+N_Q(k,-l,t)) \Big\},
\end{split}
 \end{align}
\end{subequations}
\begin{flushright}
$ -\infty<x<\infty,~0<y<\infty,~t>0.$
\end{flushright}

Denote by $H_1(k,l,t)$ and $H_2(k,l,t)$ the following functions
appearing in equations \eqref{first}:
\begin{subequations}
 \begin{align}
&H_1(k,l,t)=\f{1}{k^2+l^2} \Big\{2(\la k^3+k l^2(\la+2\mu))V^{(1)}(k,l,t)-4\mu k l^2 V^{(2)}(k,l,t)\Big\},\\
&H_2(k,l,t)=\f{1}{k^2+l^2}\Big\{4\mu k l^2 U^{(1)}(k,l,t)+2(\mu k^3-\mu k l^2) U^{(2)}(k,l,t)\Big\}.
 \end{align}
\end{subequations}
We observe that for any fixed $k\in\mathbb{R}$ and fixed $t$, $0\leq
t<T$, $T>0$, the above functions are analytic in the entire complex
$l-$plane. Indeed, equation \eqref{expansion} and the analogous
equation for $U^{(2)}$, imply that in the neighbourhood of $l=\pm
ik$, the following expansion is valid:
\begin{equation}
H_2(k,l,t)=2i\mu k \int_0^t\q{u}(k,s)ds+o(\om_1^2).
\end{equation}
Similarly,
\begin{equation}
H_1(k,l,t)=2i\la k \int_0^t\q{v}(k,s)ds+o(\om_1^2).
\end{equation}

The restriction $y>0$, as well as the analyticity of the functions
$H_1(k,l,t)$ and $H_2(k,l,t)$, allow us to deform the contour of
integration from the real axis to a contour $\ga_k$ in the upper
half $l-$plane (the particular choice of $\ga_k$ will be determined
in the next section):
\begin{subequations}\label{before_final}
\begin{align}
\label{first_uf}
\begin{split}
&u(x,y,t)=\\
&\int_{-\infty}^{\infty}\!\!dk \int_{\gamma_k} dl~\f{e^{ikx+ily}}{k^2+l^2} \Big\{2(\la k^3+k l^2(\la+2\mu) )V^{(1)}(k,l,t)-4\mu k l^2 V^{(2)}(k,l,t)\Big\}+\\
&\int_{-\infty}^{\infty}dk \int_{-\infty}^{\infty}
dl~\f{e^{ikx+ily}}{k^2+l^2}\Big\{k[N_P(k,l,t)+N_P(k,-l,t)]+\\
&~~~~~~~~~~~~~~~~~~~~~~~~~~~~~~~~~~~~~~~l[N_Q(k,l,t)-N_Q(k,-l,t)]\Big\},
\end{split}\\
\label{first_vf}
\begin{split}
&v(x,y,t)=\\
&\int_{-\infty}^{\infty}\!\!dk \int_{\gamma_k} dl~\f{e^{ikx+ily}}{k^2+l^2}\Big\{2\mu k l^2 U^{(1)}(k,l,t)+(\mu k^3-\mu k l^2) U^{(2)}(k,l,t)\Big\}+\\
&\int_{-\infty}^{\infty}dk \int_{-\infty}^{\infty}
dl~\f{e^{ikx+ily}}{k^2+l^2}\Big\{l[N_P(k,l,t)-N_P(k,-l,t)]-\\
& ~~~~~~~~~~~~~~~~~~~~~~~~~~~~~~~~~~~~~~k[N_Q(k,l,t)+N_Q(k,-l,t)]\Big\}.
\end{split}
\end{align}
\end{subequations}

\section{The Elimination of the Transforms of the Unknown Boundary Values}
Let
\begin{equation}\label{l21}
l_{21}=-l\Big(\f{\la+2\mu}{\mu}+\f{\la+\mu}{\mu}\f{k^2}{l^2}\Big)^{\f{1}{2}}
\end{equation} and
\begin{equation}\label{l12}
l_{12}=-l\Big(\f{\mu}{\la+2\mu}-\f{\la+\mu}{\la+2\mu}\f{k^2}{l^2}\Big)^{\f{1}{2}}.
\end{equation}

\begin{figure}[ht]
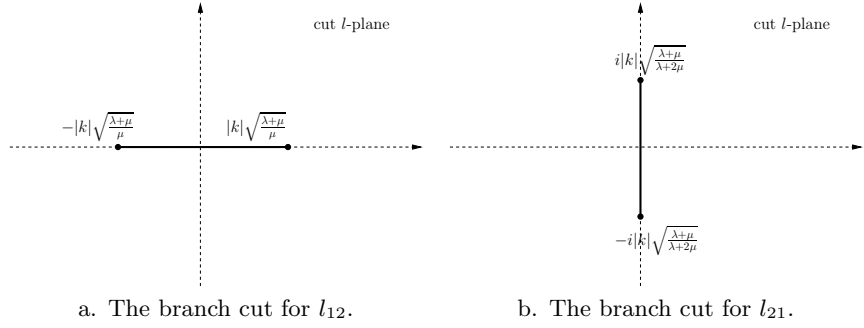

\begin{center} \small
$$
\begin{array}{cc}
\resizebox{5.5cm}{!}{\input{l12.pstex_t}} & \resizebox{5.5cm}{!}{\input{l21.pstex_t}}\\
\mbox{a. The branch cut for }l_{12}. & \mbox{b. The branch cut for }
l_{21}.
\end{array}
$$
\end{center}
\caption{Branch cuts for $l_{12}$ and $l_{21}$ in complex
$l-$plane.} \label{l12l21}
\end{figure}
The function $l_{12}$ has two branch points $\pm
k\sqrt{\f{\la+\mu}{\mu}}$, which we connect by a horizontal branch
cut; the function $l_{21}$ has two branch points $\pm
ik\sqrt{\f{\la+\mu}{\la+2\mu}}$, which we connect by a vertical
branch cut, see Fig.\ref{l12l21}. We fix the branches by the
requirements that,
$$l_{21}\sim -l\sqrt{\f{\la+2\mu}{\mu}},~~l_{12}\sim-l\sqrt{\f{\mu}{\la+2\mu}},~~ as~l\rightarrow\infty.$$

The following transformations are valid in the cut $l$-plane:
\begin{subequations}\label{superb_2}
\begin{align}
&l\rightarrow l_{21}~~ maps~~~ \om_2~~ to~~ -\om_1,~~~ U^{(2)}~~
to~~ U^{(1)},~~~ V^{(2)}~~ to~~V^{(1)};\\
&l\rightarrow l_{12}~~ maps~~~ \om_1~~ to~~ -\om_2,~~~ U^{(1)}~~
to~~ U^{(2)},~~~ V^{(1)}~~ to~~ V^{(2)}.
\end{align}
\end{subequations}

Using in equations \eqref{1} and \eqref{2} the transformations
$l\rightarrow l_{12}$ and $l\rightarrow-l$ respectively, and then
combining the two resulting equations, we obtain the following equations:
\begin{subequations}\label{C1}
 \begin{align}
 \begin{split}
k\hat{u}(k,l_{12},t)+l_{12}\hat{v}(k,l_{12},t)&=(\la k^2+l_{12}^2(\la+2\mu))V^{(2)}(k,l,t)\\
& +2\mu k l_{12}U^{(2)}(k,l,t)+N_P(k,l_{12},t),
\end{split}\\
 \begin{split}
-l\hat{u}(k,-l,t)-k\hat{v}(k,-l,t)&=-\mu(k^2-l^2)U^{(2)}(k,l,t)+\\
& 2\mu k l V^{(2)}(k,l,t) +N_Q(k,-l,t).
\end{split}
 \end{align}
\end{subequations}

Similarly, using in equations \eqref{1} and \eqref{2} the
transformations $l\rightarrow -l$ and $l\rightarrow l_{21}$
respectively, and then combining the two resulting equations, we
obtain the following equations:
\begin{subequations}\label{C2}
 \begin{align}
 \begin{split}
k\hat{u}(k,-l,t)-l\hat{v}(k,-l,t)&=(\la k^2+l^2(\la+2\mu))
V^{(1)}(k,l,t)\\
&-2\mu k lU^{(1)}(k,l,t)+N_P(k,-l,t),
\end{split}\\
\begin{split}
l_{21}\hat{u}(k,l_{21},t)-k\hat{v}(k,l_{21},t)&=-\mu(k^2-
l_{21}^2)U^{(1)}(k,l,t)\\
&-2\mu k l_{21}V^{(1)}(k,l,t)+N_Q(k,l_{21},t).
\end{split}
 \end{align}
\end{subequations}
Let
\begin{subequations}\label{Cs}
\begin{align}
&C_1(k,l)=\la k^2+l_{12}^2 (\la+2\mu),~C_2(k,l)=2\mu k
l_{12},\\
&C_3(k,l)=2\mu kl,~C_4(k,l)=-\mu(k^2-l^2),\\
&D_1(k,l)=\la k^2+l^2 (\la+2\mu),~D_2(k,l)=-2\mu kl,\\
&D_3(k,l)=-2\mu k l_{21},~D_4(k,l)=-\mu(k^2- l_{21}^2),\\
&\Delta_1(k,l)=C_1(k,l)C_4(k,l)-C_2(k,l)C_3(k,l),\\
&\Delta_2(k,l)=D_1(k,l)D_4(k,l)-D_2(k,l)D_3(k,l).
\end{align}
\end{subequations}
Simplifying the expressions for $\Delta_1$ and $\Delta_2$ we find
\begin{subequations}\label{Ts}
\begin{align}
\Delta_1(k,l)&=\mu^2(k^2-l^2)^2-4\mu^2 k^2 ll_{12},\\
\Delta_2(k,l)&=(\la k^2+l^2(\la+2\mu))^2-4\mu^2 k^2l l_{21}.
\end{align}
\end{subequations}
Equations \eqref{C1} imply
\begin{subequations}\label{X}
\begin{align}
\begin{split}
V^{(2)}(k,l,t)&=\f{1}{\Delta_1(k,l)}[C_4(k,l)(k\hat{u}(k,l_{12},t)+l_{12}\hat{v}(k,l_{12},t)-N_P(k,l_{12},t))\\
&-C_2(k,l)(-l\hat{u}(k,-l,t)-k\hat{v}(k,-l,t)-N_Q(-l))],
\end{split}\\
\begin{split}
U^{(2)}(k,l,t)&=\f{1}{\Delta_1(k,l)}[-C_3(k,l)(k\hat{u}(k,l_{12},t)+l_{12}\hat{v}(k,l_{12},t)-N_P(k,l_{12},t))\\
&+C_1(k,l)(-l\hat{u}(k,-l,t)-k\hat{v}(k,-l,t)-N_Q(k,-l,t))];
\end{split}
\end{align}
\end{subequations}
Equations \eqref{C2} imply
\begin{subequations}\label{Y}
\begin{align}
\begin{split}
V^{(1)}(k,l,t)&=\f{1}{\Delta_2(k,l)}[D_4(k,l)(k\hat{u}(k,-l,t)-l\hat{v}(k,-l,t)-N_P(k,-l,t))\\
&-D_2(k,l)
(l_{21}\hat{u}(k,l_{21},t)-k\hat{v}(k,l_{21},t)-N_Q(l_{21}))],
\end{split}\\
\begin{split}
U^{(1)}(k,l,t)&=\f{1}{\Delta_2(k,l)}[-D_3(k,l)(k\hat{u}(k,-l,t)-l\hat{v}(k,-l,t)-N_P(k,-l,t))\\
&+D_1(k,l)(l_{21}\hat{u}(k,l_{21},t)-k\hat{v}(k,l_{21},t)-N_Q(k,l_{21},t))].
\end{split}
\end{align}
\end{subequations}

\begin{figure}[ht]
\begin{center}
\resizebox{5.5cm}{!}{\input{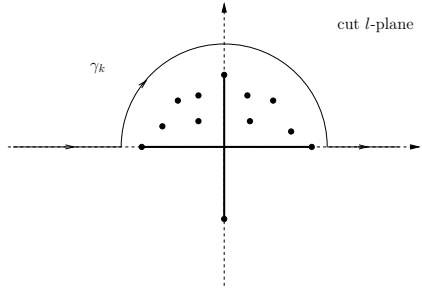}}
\end{center}
\caption{$\gamma_k$, the deformed path of integration.} \label{gamma_k}
\end{figure}

We fix the choice of the contour $\ga_k$ by requiring that every
term in the RHS of \eqref{X} and \eqref{Y} does not have a pole or a
branch point above this contour, see Fig.\ref{gamma_k}.

Regarding the zeros of $\Delta_j$, $j=1,2$, we note that they are of the form
$l=\al k$, for some constant $\al\in \mathbb{C}$. For example, if $\la=2\mu$, we find
that the zeros of $\Delta_1$ are
$$l\approx (-1.624\pm 0.126 i)k,~l\approx \pm 0.357 ik,~l\approx \pm 1.056 ik,~l \approx (1.624\pm 0.126i)k;$$
whereas the zeros of $\Delta_2$ are
$$l\approx (-0.295\pm 0.442 i)k,~l\approx \pm 0.885 ik,~l\approx \pm ik,~l \approx (0.295\pm 0.442i)k.$$

Substituting equations \eqref{X} and \eqref{Y} in equations
\eqref{before_final},
and using Jordan's lemma in the complex $l-$plane above the contour
$\ga_k$, it follows that $\hat{u}(k,-l,t)$, $\hat{u}(k,l_{12},t)$,
$\hat{u}(k,l_{21},t)$, $\hat{v}(k,-l,t)$, $\hat{v}(k,l_{12},t)$,
$\hat{v}(k,l_{21},t)$, yield a zero contribution. Hence equations \eqref{before_final} become the following equations:
\begin{subequations}\label{final}
 \begin{align}
\label{final_u}
\begin{split}
&u(x,y,t)=\\
&\int_{-\infty}^{\infty}\!\!dk \int_{\gamma_k}
\f{dl~e^{ikx+ily}}{(k^2+l^2)\Delta_1(k,l)}\Big\{2(\la k^3+k l^2(\la+2\mu))[-D_4(k,l)N_P(k,-l,t)\\
&+D_2(k,l) N_Q(k,l_{21},t)]-4\mu k l^2 [-C_4(k,l)N_P(k,l_{12},t)+C_2(k,l)N_Q(k,-l,t)]\Big\}\\
&+\int_{-\infty}^{\infty}dk \int_{-\infty}^{\infty}
dl~\f{e^{ikx+ily}}{k^2+l^2}\Big\{k[N_P(k,l,t)+N_P(k,-l,t)]+\\
&~~~~~~~~~~~~~~~~~~~~~~~~~~~~~~~~~~~~~l[N_Q(k,l,t)-N_Q(k,-l,t)]\Big\},
\end{split}\\
\label{final_v}
\begin{split}
&v(x,y,t)=\\
&\int_{-\infty}^{\infty}\!\!dk \int_{\gamma_k}
\f{dl~e^{ikx+ily}}{(k^2+l^2)\Delta_2(k,l)}\Big\{2\mu k l^2 [D_3(k,l) N_P(k,-l,t)-D_1(k,l) N_Q(k,l_{21},t)]\\
&+(\mu k^3-\mu k l^2)[C_3(k,l) N_P(k,l_{12},t)-C_1(k,l) N_Q(k,-l,t)]\Big\}\\
&+\int_{-\infty}^{\infty}dk \int_{-\infty}^{\infty}
dl~\f{e^{ikx+ily}}{k^2+l^2}\Big\{l[N_P(k,l,t)-N_P(k,-l,t)]-\\
&~~~~~~~~~~~~~~~~~~~~~~~~~~~~~~~~~~~~~k[N_Q(k,l,t)+N_Q(k,-l,t)] \Big\},
\end{split}
 \end{align}
\end{subequations}
where $\{C_j,D_j\}_1^4$ and $\{\Delta_j\}_1^2$ are defined in \eqref{Cs} and
\eqref{Ts}, and the known functions $N_P$ and $N_Q$ are defined in
\eqref{NPNQ}.

We summarize the above result in the following proposition:
\begin{proposition}
Let $(u,v)$ satisfy the Lam\'e-Navier equations \eqref{PLN} in the half plane
$$-\infty<x<\infty,~0<y<\infty,~t>0,$$
with the initial conditions \eqref{ic} and the stress boundary conditions \eqref{bv}. Assume that the given
functions $$\{u_j(x,y),v_j(x,y)\}_{j=0,1},~~\{g_j(x,t)\}_{j=1,2}$$ have sufficient smoothness and decay.
A solution of the above initial-boundary value problem, which decays for large $(x,y)$, is given by equations \eqref{final}, where:
\begin{enumerate}
\item[a)] The known functions $(N_P,N_Q)$ are defined in \eqref{NPNQ} in terms of the transforms of the initial and boundary data (see equations \eqref{dispersion}-\eqref{BigG}).
\item[b)] $\{C_j,D_j\}_1^4$ and $\{\Delta_j\}_1^2$ are given by equations \eqref{Cs} and \eqref{Ts}.
\item[c)] The contours $\gamma_{k}$ depicted in Fig.\ref{gamma_k}, are deformations of the real axis and determined by the requirement that the zeros of $\Delta_j,j=1,2$ and the branch points of $l_{21}$ and $l_{12}$ are all below $\gamma_k$.
\end{enumerate}
\end{proposition}
\section{The Normal Line Load with the Homogeneous Initial Conditions}
Consider a normal line load suddenly applied to an isotropic elastic half plane body.
In this case,
\begin{equation}\label{z_ini_normal_line_load}
u_0=u_1=v_0=v_1=0; ~g_1=0,~g_2=\sigma_0\delta(x)h(t)/(\la+\mu),
\end{equation} where $\sigma_0$ is a constant and $h(t)$ is the Heaviside
function defined by $h(t)=0,~t\leq 0;$ $h(t)=1,~t>0$.

In what follows, we compute the functions needed in equations
\eqref{final}. Equations \eqref{gtilde} and \eqref{BigG} imply
\begin{subequations}
\begin{align}
&\q{g}_1(k,t)=0,~\q{g}_2(k,t)=\f{\sigma_0}{\la+\mu}h(t),\\
&g^{(j)\pm}(k,l,t)=0,~G^{(j)}(k,l,t)=0,\\
&f^{(j)\pm}(k,l,t)=\f{\sigma_0}{\la+\mu} \f{1}{\pm i\om_j}(e^{\pm i\om_j t}-1),\\
&F^{(j)}(k,l,t)=-\f{i\sigma_0}{\la+\mu}\Big(\f{1}{\om_j^2}
-\f{\cos{(\om_j t)}}{\om_j^2}\Big), ~~~~~j=1,2. \end{align} \end{subequations}
Thus, the known functions $N_P$ and $N_Q$ defined in \eqref{NPNQ}, are given by
\begin{subequations}
\begin{align}
&N_P(k,l,t)=-l\sigma_0
\f{\la+2\mu}{\la+\mu}\Big(\f{1}{\om_1^2}-\f{\cos{(\om_1
t)}}{\om_1^2}\Big),\\
&N_Q(k,l,t)=k\sigma_0\f{\la+2\mu}{\la+\mu}\Big(\f{1}{\om_2^2}-\f{\cos{(\om_2
t)}}{\om_2^2}\Big).
\end{align}
\end{subequations}
The above equations imply
\begin{subequations}
\begin{align}
&N_P(k,l,t)+N_P(k,-l,t)=0,\\
&N_Q(k,l,t)+N_Q(k,-l,t)=2k\sigma_0\f{\la+2\mu}{\la+\mu}\Big(\f{1}{\om_2^2}-\f{\cos{(\om_2 t)}}{\om_2^2}\Big),\\
&N_P(k,l,t)-N_P(k,-l,t)=-2l\sigma_0
\f{\la+2\mu}{\la+\mu}\Big(\f{1}{\om_1^2}-\f{\cos{(\om_1 t)}}{\om_1^2}\Big),\\
&N_Q(k,l,t)-N_Q(k,-l,t)=0,\\
&N_P(k,l_{12},t)=-\f{l_{12}}{k} N_Q(k,l,t),\\
&N_Q(k,l_{21},t)=-\f{k}{l}N_P(k,l,t).
\end{align}
\end{subequations}
Substituting the above relations and equations \eqref{Cs} and \eqref{Ts} into equations
\eqref{final}, we find the following result:
\begin{proposition}
Let $(u,v)$ satisfy the Lam\'e-Navier equations \eqref{PLN} in the half plane
$$-\infty<x<\infty,~0<y<\infty,~t>0,$$
with homogeneous initial condition and the normal line load boundary condition \eqref{z_ini_normal_line_load}. The solution of this initial-boundary value problem, which decays for large $(x,y)$, is given by
\begin{subequations}
\begin{align}
\begin{split}
& u(x,y,t)=\\
&-\sigma_0 \f{\la+2\mu}{\la+\mu}\int_{-\infty}^{\infty}dk
\int_{\gamma_k}
dl~\f{e^{ikx+ily}}{(k^2+l^2)[\mu^2(k^2-l^2)^2-4\mu^2 k^2 ll_{12}]}\\
&\Big\{2(\la k^3l+k l^3(\la+2\mu))(1-\cos{(\om_1 t)})+4\mu k l^2
l_{12}(1-\cos{(\om_2 t)})\Big\},
\end{split}\\
\begin{split}
&v(x,y,t) =\\
&-\sigma_0 \f{\la+2\mu}{\la+\mu}\int_{-\infty}^{\infty}dk
\int_{\gamma_k}
dl~\f{e^{ikx+ily}}{(k^2+l^2)[(\la k^2+l^2(\la+2\mu))^2-4\mu^2 k^2l l_{21}]}\\
&\Big\{2\mu k^2 l^2[2\mu ll_{21}+(\la k^2+l^2(\la+2\mu))]\f{1-\cos{(\om_1 t)}}{\om_1^2}\\
&~~~~~~~~~~~~~~~~~~~~~~~~~~~~~-2\mu^2(k^4-k^2 l^2)(k^2-l^2-2ll_{12})\f{1-\cos{(\om_2 t)}}{\om_2^2}\Big\}\\
&-2\sigma_0\f{\la+2\mu}{\la+\mu}\int_{-\infty}^{\infty}dk
\int_{-\infty}^{\infty}
dl~\f{e^{ikx+ily}}{k^2+l^2}\Big\{l^2\f{1-\cos{(\om_1
t)}}{\om_1^2}+k^2\f{1-\cos{(\om_2 t)}}{\om_2^2}\Big\},
\end{split}
\end{align}
\end{subequations}
where $(\om_1,\om_2)$ are defined by equations \eqref{dispersion_choice} and $l_{12},l_{21}$ are defined by equations \eqref{l12}, \eqref{l21} respectively.
\end{proposition}

Similar results can be obtained for the other Lamb's problems.

\section{Conclusions}
The main result of this paper is the derivation of equations
\eqref{final}. These equations express the displacements
$(u(x,y,t),v(x,y,t))$ in terms of an integral along the real line
and integrals along contours $\ga_k$ of the complex $l-$plane; these
integrals involve transforms of the given initial and boundary
data. Indeed, equations \eqref{final} involve the functions
$N_P(k,l,t)$ and $N_Q(k,l,t)$, which are defined in equations
\eqref{NPNQ} in terms of the Fourier transforms $P_0(k,l,t),~
Q_0(k,l,t),$ $P_1(k,l,t),~ Q_1(k,l,t)$ of the initial data
(see equations \eqref{initial}), as well as in terms of certain
known transforms $G^{(1)}(k,l,t),$ $G^{(2)}(k,l,t)$,
$F^{(1)}(k,l,t),$ $F^{(2)}(k,l,t)$ of the boundary data (see
equations \eqref{BigG}).

The starting point of the derivations of equations \eqref{final}
is the derivation of the global relations \eqref{global_relations}. These equations are the direct
consequence of the application of the two-dimensional Fourier transform
and of the substitution in the resulting equations of the given initial
and boundary conditions. Equations \eqref{global_relations} involve the \textit{known}
functions $N_P$ and $N_Q$, as well as the \textit{unknown} functions
$U^{(1)},U^{(2)},V^{(1)},V^{(2)}$ (these functions involve certain
transforms of the unknown boundary values $u(x,0,t)$ and $v(x,0,t)$, see
equations \eqref{xF} and \eqref{BigU}). The elimination of the above unknown
functions is achieved by the following steps:
1. By employing the transformation $l\rightarrow -l$, which leave the unknown
functions $(U^{(1)},U^{(2)},V^{(1)},V^{(2)})$ invariant, and by utilising the analyticity
properties of these functions, we obtain the integral representations
\eqref{before_final}. These representations involve an integral along the real $k-axis$
and an integral along the contour $\gamma_k$ of the complex $l-$plane.
2. Using the transformations $l\rightarrow l_{21}$ and $l\rightarrow l_{12}$
which map $U^{(2)}$ to $U^{(1)}$,$V^{(2)}$ to $V^{(1)}$ and
$U^{(1)}$ to $U^{(2)}$,$V^{(1)}$ to $V^{(2)}$ respectively, we
express the unknown functions $(U^{(1)},U^{(2)},V^{(1)},V^{(2)})$ in
terms of known functions as well as in terms of certain unknown
functions which however are analytic in a certain domain of the complex $l-$plane,
see equations \eqref{X} and \eqref{Y}. 3. Using \eqref{X} and \eqref{Y} in equations \eqref{before_final} and
employing Jordan's lemma we obtain equations \eqref{final}.

The main advantages of the new approach are the following:
\begin{enumerate}
\item{The new method provides an analytic solution of Lamb's problem with
arbitrary initial and boundary conditions.}
\item{This solution is expressed in terms of the given initial and
boundary data. The relevant representation is novel even for the particular case of
homogeneous initial conditions (this case has been analysed by several
authors). An alternative approach using the Laplace transform is briefly discussed in the Appendix. By comparing equations \eqref{final} and equation \eqref{int2}, the advantage of the new formulae becomes clear. Actually, taking into consideration that the initial-boundary value problems of the equations of elastodynamics are well posed for any finite $t$, the use of the Laplace transform, which requires $t\rightarrow\infty$, is clearly inappropriate.}
\item{The new method can be employed for the solution of several related initial-boundary value problems,
including the problem of the isotropic half space in the axisymmetric case.
Furthermore, it can be extended to problems in three dimensions \cite{Fokas_Yang_3D}.}
\end{enumerate}

\section*{Appendix}
It is possible to analyse Lamb's problem by using only the transformations \eqref{superb_1} instead of using the transformations \eqref{superb_1} \textit{and} the transformations \eqref{superb_2}. However, in this case one \textit{cannot} eliminate directly all unknown boundary values. Instead, one can derive a complicated expression for the unknown boundary values in terms of the given initial and boundary data. (This approach is similar with the one used in \cite{Ashton_ASF} for solving Crighton's problem).
\begin{figure}[ht]
\begin{center}
\resizebox{5.5cm}{!}{\input{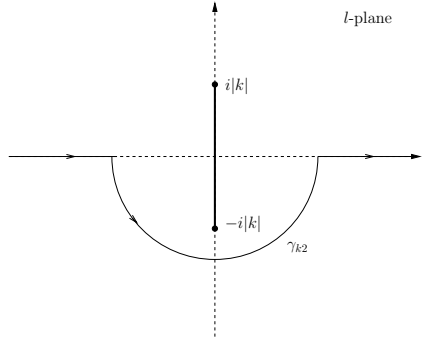}}
\end{center}
\caption{$\gamma_{k2}$, the path of integration.} \label{fig:path2}
\end{figure}
Indeed, let the contour $\gamma_{k2}$ be a simple curve in the lower
half $l-$plane, determined by the requirement that it does not cross
the branch cut associated with $\om_1$ and $\om_2$, see
Fig.\ref{fig:path2}. Let $K$ be the integral operator defined by
$$(K[f])(k,t)=\f{1}{2\pi}\int_{\gamma_{k2}}\f{f(k,l,t)~dl
}{l(k^2+l^2)^{1/2}},$$
for any function $f(k,l,t)$ with appropriate smoothness and decay. Integrating the global relations
\eqref{global_relations} along $\gamma_{k2}$ we find that the functions
$$\q{h}(k,t)=\left(
\begin{array}{llcl}
\q{u}(k,t)\\
\q{v}(k,t)\\
\end{array}\right),~~~
\q{g}(k,t)=\left(
\begin{array}{llcl}
\q{g}_1(k,t)\\
\q{g}_2(k,t)\\
\end{array}\right),
$$
satisfy a system of
Volterra integral equations of the second kind:
\begin{equation}\label{int2}
\q{h}(k,t)=(K[N])(k,t)\ast\q{g}(k,t)+(K[M])(k,t)\ast \q{h}(k,t)+(K[H])(k,t),
\end{equation}
\begin{flushright}
$k\in\mathbb{R}, 0\leq t<T,$
\end{flushright}
where $T>0$, $\ast$ denotes the convolution operation with respect to
$t$,
$$N(k,l,t)=\left(
\begin{array}{llcl}
i\sqrt{\mu}le^{i\om_2t}                & -i\sqrt{\mu}\f{\la+2\mu}{\mu}k e^{i\om_2t}\\
i\f{\mu}{\sqrt{\la+2\mu}}ke^{i\om_1t} & i\sqrt{\la+2\mu}le^{i\om_1t}\\
\end{array}\right),$$
$$M(k,l,t)=\left(
\begin{array}{llcl}
\sqrt{\mu}k^2e^{i\om_2t}                & 2\sqrt{\mu}kl e^{i\om_2t}\\
-\f{2\mu}{\sqrt{\la+2\mu}}kle^{i\om_1t} & -\f{\la}{\sqrt{\la+2\mu}}k^2 e^{i\om_1t}\\
\end{array}\right),$$ and the known function $H(k,l,t)$ is defined by
$$H(k,l,t)=\Big(i\f{e^{i\om_1t}}{\la+2\mu}(i\om_1P_0(k,l)+P_1(k,l)),~i\f{e^{i\om_2t}}{\mu}(i\om_2Q_0(k,l)+Q_1(k,l))\Big)^T.$$
The solution of the integral equations
\eqref{int2} yields the unknown transforms appearing in
\eqref{global_relations}.

Equations \eqref{int2} can be solved in closed form by using the
Laplace transform in $t$. Let ``$\circ$" denote the Laplace
transform with respect to $t$ and let
\begin{equation}
w_1=(p^2+(\la+2\mu)k^2)^{\f{1}{2}},~w_2=(p^2+\mu k^2)^{\f{1}{2}},
\end{equation}
\begin{equation}
\Delta=\f{(p^2+2\mu
k^2)^2}{w_1^2w_2^2}-4\sqrt{\f{\mu}{\la+2\mu}}\f{\mu k^2}{w_1 w_2}.
\end{equation}
The solution of \eqref{int2} is given by
\begin{equation}\label{Spectral_NDC}
\begin{split}
&\mathring{\q{h}}=\left(\begin{array}{llcl}
-\f{\sqrt{\mu}p^2 w_2}{(p^2+2\mu k^2)^2-4\sqrt{\f{\mu}{\la+2\mu}}\mu k^2w_1w_2}  & ik\f{(\la+2\mu)(p^2+2\mu k^2)-2\sqrt{\mu(\la+2\mu)}w_1w_2}{(p^2+2\mu k^2)^2-4\sqrt{\f{\mu}{\la+2\mu}}\mu k^2w_1w_2}\\
ik \f{2\mu\sqrt{\f{\mu}{\la+2\mu}}w_1w_2-\mu(p^2+2\mu k^2)}{(p^2+2\mu k^2)^2-4\sqrt{\f{\mu}{\la+2\mu}}\mu k^2w_1w_2} & -\f{\sqrt{\la+2\mu} p^2 w_1}{(p^2+2\mu k^2)^2-4\sqrt{\f{\mu}{\la+2\mu}}\mu k^2w_1w_2}\\
\end{array}\right)\mathring{\q{g}}\\
&+\f{1}{\Delta} \left(\begin{array}{llcl}
\f{p^2+2\mu k^2}{w_1^2}               & 2i\f{\sqrt{\mu}k}{w_2}\\
-2i\sqrt{\f{\mu}{\la+2\mu}}\f{\sqrt{\mu}k}{w_1} & \f{p^2+2\mu k^2}{w_2^2}\\
\end{array}\right)
\left(\begin{array}{llcl}
&-\f{i}{2\pi}\int_{\gamma_{k2}}\f{i\om_1P_0+P_1}{(\la+2\mu)(p-i\om_1)}\f{dl}{l(k^2+l^2)^{\f{1}{2}}}\\
&\f{i}{2\pi}\int_{\gamma_{k2}}\f{i\om_2Q_0+Q_1}{\mu(p-i\om_2)}\f{dl}{l(k^2+l^2)^{\f{1}{2}}}\\
\end{array}\right).
\end{split}
\end{equation}
The zeros of $\Delta$ coincide with the zeros of Rayleigh's
function. When $\mu/\la>0.906$, the known
transforms $\mathring{\q{g}}(k,p)$, for each fixed $k$, do \emph{not} have poles with
positive real parts.

In the particular case of the problem of homogeneous initial condition and the normal line
load boundary condition, equation \eqref{Spectral_NDC} reduces to
the classic Lamb's solution \cite{HL}\cite{JM}.

\paragraph{Acknowledgements}
We would like to thank Ishan Sharma for drawing our attention to
Lamb's problems. We are grateful to Michail Dimakos, Dionyssis
Mantzavinos and Anthony Ashton for several valuable observations. Di
Yang would like to thank Professor Youjin Zhang for his advise and
for several discussions, as well as the China Scholarship Council
for supporting him for a joint PhD study at the University of
Cambridge. A. S. Fokas would like to thank Y. Antipov for his
suggestions and also acknowledges the generous support of the
Guggenheim Foundation, USA. The work of Yang is partially supported
by the National Basic Research Program of China (973 Program)
No.2007CB814800.

\bibliographystyle{amsplain}

\begin{thebibliography}{10}
\bibitem{HL}
H. Lamb, {\em On the propagation of tremors over the surface of an
Elastic Solid}, Phil. Trans. Royal. Soc. London, Series A, Vol. 203, pp.
1-42, 1904.

\bibitem{JM}
J. Miklowitz, {\em The theory of elastic waves and waveguides},
North-Holland Publishing Company, 1978.

\bibitem{KLJ}
K. L. Johnson, {\em Contact mechanics,} Cambridge University Press,
Cambridge, 1985.

\bibitem{point_load0}
C. L. Pekeris, {\em The seismic surface pulse,} Proc. Nat. Acad.
Sci., vol. 41, 1955.

\bibitem{point_load05}
R. G. Payton, {\em An application of the dynamic Betti-Rayleigh
reciprocal theorem to moving-point loads in elastic media,} Quart.
Appl. Math. Vol. 21, pp. 299-313, 1964.

\bibitem{point_load1}
K. Schiel and S. J. Prot\'azio, {\em Transient-wave solution for
Lamb's problem at the free surface}, Bulletin of the Seismological
Society of America, Vol. 79, No. 6, pp. 1956-1971, Dec., 1989.

\bibitem{point_load2}
D. C. Gakenheimer, {\em Transient excitation of an elastic-half space
by a point load travelling on the surface,} PhD Thesis, California
Institute of Technology, 1969.

\bibitem{point_load4}
J. R. Barber, {\em Surface displacements due to a steadily moving
point force,} J. Appl. Mech., Vol. 63, Jun., 1976.


\bibitem{point_load5}
M. C. M. Bakkeraq, M. D. Verweij, B. J. Kooij, H. A. Dietermana,
{\em The travelling point load revisited}, Wave Motion 29, 1999.

\bibitem{line_load1}
R. G. Payton, {\em Transient motion of an elastic half-space due to
a moving surface line loads,} Int. J. Engng Sci. Vol. 5, pp. 49-79,
1967.


\bibitem{line_load16}
J. Miklowitz, W. R. Garrott, {\em Lamb's problem for an impulsive
line load on a Laminated Composite}, Modern problems in elastic wave
propagation, IUTAM, Symposium at Northwestern University, Illinois,
USA, 1977.

\bibitem{line_load2}
H. G. Georgiadis and J. R. Barber, {\em Steady-state transonic
motion of a line load over an elastic half-space. The corrected
Cole/Huth solution,} Journal of Applied Mechanics, Vol. 60, No. 3,
772-774, 1993.

\bibitem{Payton}
R. G. Payton, {\em Elastic wave propagation in transversely isotropic
media}, Martinus Nijhoff Publishers, USA, 1983.

\bibitem{ASF1997}
A. S. Fokas, {\em A unified transform method for solving linear and
certain nonlinear PDEs,} Proc. Math. Phys. Eng. Sci. 453, pp.
1411-1443, 1997.

\bibitem{ASF}
A. S. Fokas, {\em A unified approach to boundary value problems,}
SIAM, 2008.

\bibitem{Crowdy_ASF}
D. G. Crowdy, A. S. Fokas, {\em Explicit integral solutions for the plane elastostatic semi-strip}, Proc. Roy. Soc. London Ser. A 460, pp. 1289-1309, 2004.

\bibitem{Avr-ASF}
D. ben-Avraham, A. S. Fokas, {\em The solution of the modified Helmholtz equation in a wedge and an application to difussion-limited coalescence}, Phys. Lett. A 263, pp. 355-359, 1999.

\bibitem{Avr-ASF2}
D. ben-Avraham, A. S. Fokas, {\em The modified Helmholtz equation in a triangular domain and an application to diffusion-limited coalescence}, Phys. Rev. E 64, 2001.

\bibitem{ASF_Kapaev0}
A. S. Fokas, A. A. Kapaev, {\em A Riemann-Hilbert approach to the Laplace equation}, J. Math. Anal. and Appl. 251, pp. 770-804, 2000.

\bibitem{ASF_l_nonlinear}
A. S. Fokas, {\em On the integrability of linear and nonlinear PDEs}, J. Math. Phys. 41, pp. 4188-4237, 2000.

\bibitem{ASF_2d}
A. S. Fokas, {\em Two dimensional linear PDE's in a convex polygon}, Proc. R. Soc. Lond. A 457, pp. 371-393, 2001.

\bibitem{ASF_Zyskin}
A. S. Fokas, M. Zyskin, {\em The fundamental differential form and boudary value prolems}, Quart. J. Mech. Appl. Math. 55, pp. 457-479, 2002.

\bibitem{Fokas_trans}
A. S. Fokas, {\em A new transform method for evolution PDEs}, IMA J. Appl. Math. 67, pp. 1-32, 2002.

\bibitem{Fokas_var_coe}
A. S. Fokas, {\em Boudary-value problems for linear PDEs with variable coefficients,} Proc. R. Soc. London A 460, pp. 1131-1151, 2004.

\bibitem{Tre_ASF}
P. A. Treharne, A. S. Fokas, {\em Boundary-value problems for systems of evolution equations}, IMA J. Appl. Math 69, 2004.

\bibitem{Pelloni_ASF}
A. S. Fokas, B. Pelloni, {\em Boundary value problems for Boussinesq type systems}, Math. Phys. Anal. Geom. 8, pp. 59-96, 2005.

\bibitem{Pelloni_ASF2}
A. S. Fokas, B. Pelloni, {\em A transform method for evolution PDEs on the interval}, IMA J. Appl. Maths 75, pp. 564-587, 2005.

\bibitem{DTP_ASF}
A. S. Fokas, D. T. Papageorgiou, {\em Absolute and convective instability for evolution PDEs on the half-line}, Studies in Appl. Maths. 114, pp. 95-114, 2005.

\bibitem{ASF_DAP}
A. S. Fokas, D. A. Pinotsis, {\em The Dbar formalism for certain linear non-homogeneous elliptic PDEs in two dimensions}, Eur. J. Appl. Math. 17, I. 3, pp. 323-346, 2006.

\bibitem{ASF_DAP2}
A. S. Fokas, D. A. Pinotsis, {\em Quaternions, evaluation of integrals and boundary value problems}, Comp. Meth. Funct. Th. 7, pp. 443-476, 2007.

\bibitem{Tre_ASF2}
P. A. Treharne, A. S. Fokas, {\em Initial-boundary value problems for linear PDEs with variable coefficients}, Camb. Phil. Soc. 143, pp. 221-242, 2007.

\bibitem{Dassios_ASF_1}
G. Dassios, A.S. Fokas, {\em The basic elliptic equations in an equilateral triangle}, Proc R Soc A, 461, pp. 2721-2748,
2005.

\bibitem{Dassios_ASF_3}
G. Dassios, A.S. Fokas, {\em Methods for Solving Elliptic PDEs in Spherical Coordinates}, SIAM J. APPL. MATH. Vol. 68, No. 4, pp. 1080-1096, 2008.

\bibitem{Spence_ASF1}
E. A. Spence, A. S. Fokas, {\em A new transform method I: domain-dependent fundamental solutions and integral representations}, Proc R Soc A, 466, pp. 2259-2281, 2010.

\bibitem{Spence_ASF2}
E. A. Spence, A. S. Fokas, {\em A new transform method II: the global relation, and boundary value problems in polar co-ordinates}, Proc R Soc A, 466, pp. 2283-2307, 2010.

\bibitem{Pelloni}
B. Pelloni, {\em The spectral representation of two-point boundary value problems for linear evolution equations}, Proc R Soc London Ser A, 461, pp. 2965-2984, 2005.

\bibitem{Antipov_ASF}
Y. Antipov, A. S. Fokas, {\em A transform method for the modified Helmholtz equation on the semi-strip}, Math. Proc. Cambridge Philos. Soc. 137, pp. 339-365, 2004.

\bibitem{ASF_Kapaev}
A. S. Fokas, A. A. Kapaev, {\em On a transform method for the Laplace equation in a polygon}, IMAJ. Appl. Math. 68, pp. 355-408, 2003.

\bibitem{FFX}
S. Fulton, A. S. Fokas, C. Xenophontos, {\em An analytical method for lienar elliptic PDEs and its numerical implementation}, J. Comput. Appl. Math. 167, pp. 465-483, 2004.

\bibitem{SFFS}
A. G. Sifalakis, A. S. Fokas, S. R. Fulton, Y. G. Saridakis, {\em The generalised Dirichlet-Neumann map for linear elliptic PDE's and its numerical implementation}, J. Comput. Appl. Math. 219, pp. 9-34, 2008.

\bibitem{FF1}
N. Flyer, A. S. Fokas, {\em A hybrid analytical numerical method for solving evolution partial differential equations. I. The half-line}, Proc. R. Soc. 464, pp. 1823-1849, 2008.

\bibitem{FF2}
N. Flyer, A. S. Fokas, S. A. Smitheman, E. A. Spence, {\em A semi-analytical numerical method for solving evolution and elliptic partial differential equations}, J. Comp. Appl. Math. 227, pp. 59-74, 2009.

\bibitem{KF}
K. Kalimeris, A. S. Fokas, {\em The heat equation in the interior of an equilateral triangle}, Stud. Appl. Math. 124, pp. 283-305, 2010.

\bibitem{SSF}
S. A. Smitheman, E. A. Spence,  A. S. Fokas, {\em A spectral collocation method for the Laplace and modified Helmholtz equations in a convex polygon}, IMA J. Num. Anal. (in press), 2010.

\bibitem{Ashton_ASF}
A. C. L. Ashton, A. S. Fokas, {\em A novel method of solution for
the fluid-loaded plate}, Proc. R. Soc. A 465, pp. 3667-3685, 2009.

\bibitem{Fokas_Yang_3D}
A. S. Fokas, D. Yang, {\em A novel approach to elastodynamics: II. The three dimensional case}, (preprint), 2010.
\end{thebibliography}

\end{document}